\documentclass{amsart}
\usepackage[utf8]{inputenc}
\usepackage{yufei}
\newcommand{\ex}{\text{ex}}

\title{note on the Tur\'an number of the linear $3$-graph $C_{13}$}
\author{Chaoliang Tang}
\address{School of Mathematical Sciences, Fudan University,Shanghai, China 200433}
\email{cltang17@fudan.edu.cn}

\author{Hehui Wu}
\address{Shanghai Center for Mathematical Sciences, Fudan University,Shanghai, China 200438}
\email{hhwu@fudan.edu.cn}

\author{Shengtong Zhang}
\date{September 2021}
\address{Department of Mathematics, Massachusetts Institute of Technology, Cambridge, MA 02139}
\email{stzh1555@mit.edu}

\author{Zeyu Zheng}
\address{School of Mathematical Sciences, Fudan University, Shanghai, China 200433}
\email{zeyuzheng19@fudan.edu.cn}

\begin{document}
\begin{abstract}
    Let the crown $C_{13}$ be the linear $3$-graph on $9$ vertices $\{a,b,c,d,e,f,g,h,i\}$ with edges 
    $$E = \{\{a,b,c\}, \{a, d,e\}, \{b, f, g\}, \{c, h,i\}\}.$$ Proving a conjecture of Gy\'arf\'as et. al., we show that for any crown-free linear $3$-graph $G$ on $n$ vertices, its number of edges satisfy
    $$\lvert E(G) \rvert\leq \frac{3(n - s)}{2}$$
    where $s$ is the number of vertices in $G$ with degree at least $6$. This result, combined with previous work, essentially completes the determination of linear Tur\'an number for linear $3$-graphs with at most $4$ edges.
\end{abstract}
\maketitle
\section{Introduction}
A \textbf{linear $3$-graph} $G = (V, E)$ consists of a finite set of vertices $V = V(G)$ and a collection $E = E(G)$ of $3$-element subsets of $V$(edges), such that any two edges in $E$ share at most one vertex. If $H$ and $F$ are linear $3$-graphs, then $H$ is $F$-free if it contains no copy of $F$. For a linear $3$-graph $F$, and a positive integer $n$, the \textbf{linear Tur\'an number} $\ex(n, F)$ is the maximum number of edges in any $F$-free linear $3$-graph on $n$ vertices.

Let the \textbf{crown} $C_{13}$ be the linear $3$-graph on $9$ vertices $\{a,b,c,d,e,f,g,h,i\}$ with edges 
$$E = \{\{a,b,c\}, \{a, d,e\}, \{b, f, g\}, \{c, h,i\}\}.$$ 
\begin{figure}[h]
  \includegraphics[width=0.4\linewidth]{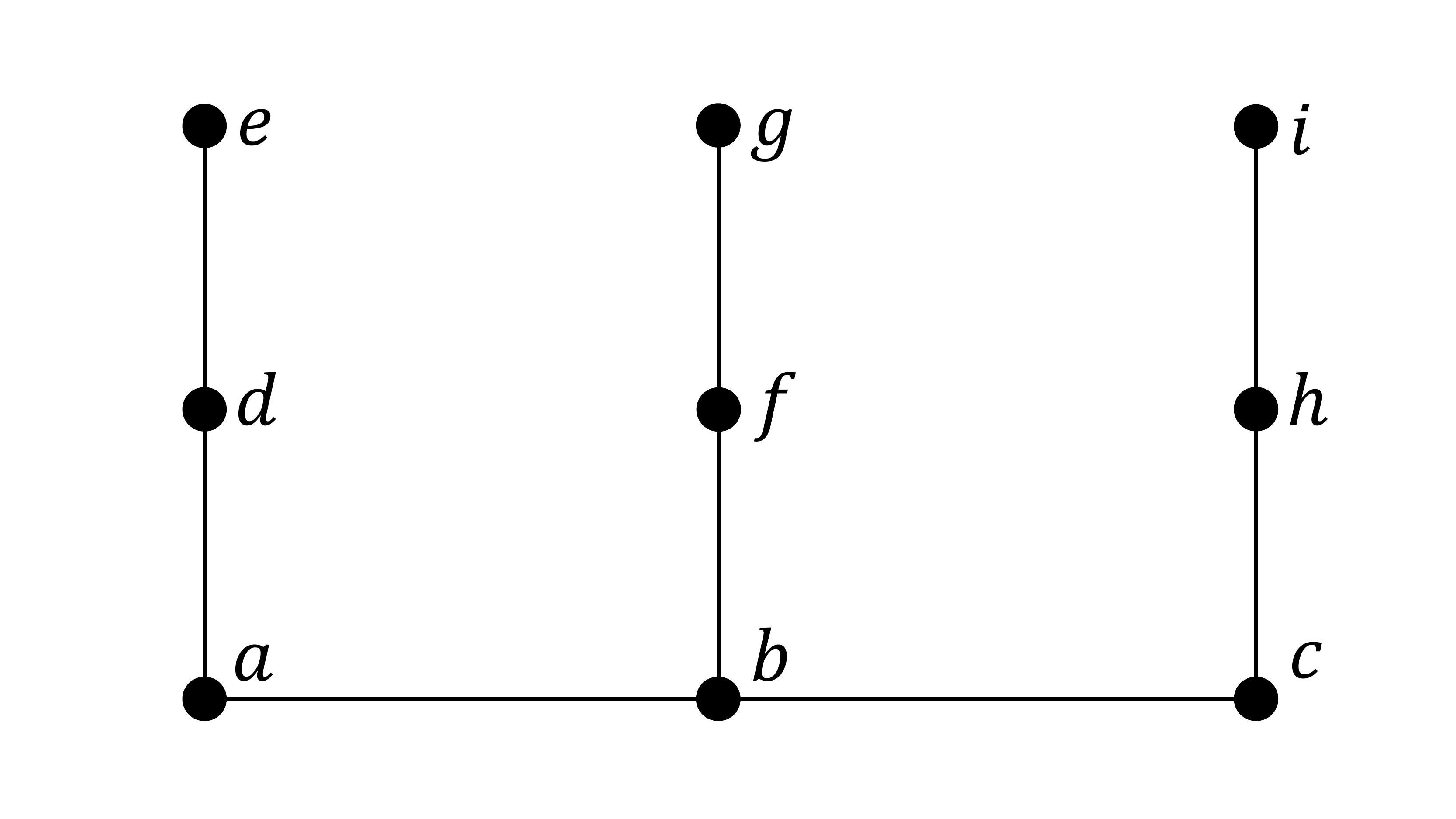}
  \caption{The crown $C_{13}$.}
\end{figure}

The study of $\ex(n, C_{13})$ was initiated by Gy\'arf\'as, Ruszink\'o and S\'ark\"ozy in \cite{G21}, where they showed the bounds
$$6\floor{\frac{n - 3}{4}} + \epsilon \leq \ex(n, C_{13}) \leq 2n.$$
where $\epsilon = 0$ if $n - 3 \equiv 0,1\bmod{4}$, $\epsilon = 1$ if $n - 3 \equiv 2\bmod{4}$, and $\epsilon = 3$ if $n - 3 \equiv 3\bmod{4}$. In \cite{C21}, Gy\'arf\'as et. al. showed that every linear $3$-graph with minimum degree $4$ contains a crown. They also proposed some ideas to obtain the exact bounds. Very recently, Fletcher showed in \cite{F21} the improved upper bound
$$\ex(n, C_{13}) < \frac{5}{3}n.$$
In this paper, we show that the lower bound in \cite{G21} is essentially tight, thus resolving a conjecture in \cite{C21}. In fact, we show the following stronger result.
\begin{theorem}
\label{thm:main1}
Let $G$ be any crown-free linear $3$-graph $G$ on $n$ vertices. Then its number of edges satisfies
$$\abs{E(G)} \leq \frac{3(n - s)}{2}.$$
where $s$ is the number of vertices in $G$ with degree at least $6$.
\end{theorem}
Furthermore, we show that when $s$ is small, the upper bound can be improved.
\begin{theorem}
\label{thm:main2}
Let $G$ be any crown-free linear $3$-graph $G$ on $n$ vertices, and let $s$ be the number of vertices in $G$ with degree at least $6$. If $s \leq 2$, then the number of edges satisfies
$$\abs{E(G)} \leq \frac{10(n - s)}{7}.$$
\end{theorem}
Combining the two theorems above, we immediately conclude that the lower bound in \cite{G21} is exact when $n \equiv 3\bmod{4}$ and $n \geq 63$.
\begin{corollary}
If $n \geq 63$, then
$$\ex(n, C_{13}) \leq \frac{3(n - 3)}{2}.$$
\end{corollary}
The paper is structured as follows. In \cref{sec:main1} and \cref{sec:main2} we present the main innovative inequality and prove our main theorems, quotient a technical and familiar lemma that we prove in \cref{sec:554}.
\section{Proof of \texorpdfstring{\cref{thm:main1}}{Theorem 1.1}}
\label{sec:main1}
Let $G$ be any linear $3$-graph. For each $v \in V(G)$, let $d(v)$ be the degree of $v$, which is the number of edges in $E(G)$ that contains $v$. For each edge $e \in E(G)$ and positive integers $a \geq b \geq c$, we write $D(e) \geq \langle a,b,c\rangle$ if we can write $e = \{x,y,z\}$ such that $d(x) \geq a$, $d(y) \geq b$ and $d(z) \geq c$.

Suppose the contrary. Let $G$ be the smallest linear $3$-graph such that $G$ has greater than $3(n - s) / 2$ edges. For each $v \in V(G)$,  let $\chi(v) = 1$ if $d(v) \leq 5$, and $\chi(v) = 0$ otherwise.

Our key innovation is the following observation
$$\sum_{e \in E(G)} \sum_{v \in V(G), v \in e} \frac{\chi(v)}{d(v)} = \sum_{v \in V(G)} \sum_{e \in E(G), v \in e} \frac{\chi(v)}{d(v)} = \sum_{v \in V(G)} \chi(v) = n - s.$$
As $\abs{E(G)} > 3(n - s)/2$, we conclude that there exists an edge $e = \{x,y,z\}$ such that
\begin{align}
\label{ineq:main}
 \frac{\chi(x)}{d(x)} + \frac{\chi(y)}{d(y)} + \frac{\chi(z)}{d(z)} < \frac{2}{3}.   
\end{align}
Without loss of generality, assume $d(x) \leq d(y) \leq d(z)$. First we note that $d(x) \geq 2$ and $d(y) \geq 4$, as otherwise \eqref{ineq:main} would be violated. If $d(z) \geq 6$, then we can easily find a $C_{13}$ by choosing an edge $e_1 \neq e$ adjacent to $x$, choosing an edge $e_2$ adjacent to $y$ that does not share a vertex with $e_1$, and finally choosing an edge $e_3$ adjacent to $z$ that does not share a vertex with $e_1$ and $e_2$, contradiction. Therefore, we have $d(z) \leq 5$, and \eqref{ineq:main} implies that $D(e) \geq \langle 5,5,4\rangle.$

We use the following lemma to handle the case $D(e) \geq \langle 5,5,4\rangle$. As the lemma is quite straightforward using the techniques in \cite{C21}, \cite{F21} and \cite{G21}, we delay the lengthy proof to \cref{sec:554}. 
\begin{restatable}{lemma}{ll}
\label{lem:554}
Let $G$ be a crown-free graph and $e = \{x,y,z\} \in E(G)$ satisfy $D(e) \geq \langle 5,5,4\rangle$. Then, the vertex set of all vertices sharing an edge with $\{x,y,z\}$,
$$S = \bigcup_{f \in E(G), f \cap \{x,y,z\} \neq \emptyset} f,$$
contains exactly $11$ vertices and all vertices in $S$ have degree at most $5$. The set of edges that contains at least one vertex in $S$,
$$E_S = \{f: f \in E(G), f \cap S \neq \emptyset\},$$
contains at most $13$ edges, and all elements of $E_S$ are subsets of $S$. In other words, the subgraph $G[S]$ is a connected component of $G$.
\end{restatable}
Let $G - S$ be the graph obtained by deleting the vertices $S$ and the edges in $E_S$. By the lemma, the graph $G - S$ has $n' = n - 11$ vertices and at least $\abs{E(G)} - 13$ edges. Furthermore, the number of vertices in $G - S$ of degree at least $6$ is exactly $s$. Therefore, we conclude that
$$\abs{E(G - S)} \geq \abs{E(G)} - 13 > \frac{3(n - s)}{2} - 13 > \frac{3(n' - s)}{2}$$
contradicting the assumption that $G$ is the smallest counterexample to \cref{thm:main1}. So we have shown \cref{thm:main1}.
\section{Proof of \texorpdfstring{\cref{thm:main2}}{Theorem 1.2}}
%% TODO
\label{sec:main2}

We use the same notations as \cref{sec:main1}.

Suppose the contrary. Let $G$ be the smallest linear 3-graph such that $G$ has at most $2$ vertices with degree at least $6$ and has greater than $10(n - s)/7$ edges. 

For each $e \in E(G)$ and $v \in e$, we define a weight $\chi(v, e)$ as follows: let  $\chi(v,e)=1$ if $d(v)=1,2,4,5$, and $\chi(v,e)=0$ if $d(v)\geq 6$. If $d(v)=3$, let $\chi (v,e)=1.05$ if there exists at least one vertex in $e$ with degree at least $6$, and $\chi(v,e)=0.9$ otherwise.

Since $s \leq 2$, we have
$$\sum_{e \in E(G)} \sum_{v \in V(G), v \in e} \frac{\chi(v,e)}{d(v)} = \sum_{v \in V(G)} \sum_{e \in E(G), v \in e} \frac{\chi(v,e)}{d(v)} \leq n - s.$$
As $\abs{E(G)} > 10(n - s)/7$, we conclude that there exists an edge $e = \{x,y,z\}$ such that
\begin{align}
\label{ineq:main2}
 \frac{\chi(x,e)}{d(x)} + \frac{\chi(y,e)}{d(y)} + \frac{\chi(z,e)}{d(z)} < \frac{7}{10}.   
\end{align}
Without loss of generality, assume $d(x) \leq d(y) \leq d(z)$. First we note that $d(x) \geq 2$, as otherwise \eqref{ineq:main2} would be violated. Then note that if $d(y) \leq 3$, no matter $d(z)$ is greater than $6$ or not \eqref{ineq:main2} would also be violated, thus $d(y) \geq 4$. 

The rest of the proof proceeds exactly the same as  \cref{sec:main1}, other than the following inequality which leads to contradiction. \cref{thm:main2} then follows.
$$\abs{E(G - S)} \geq \abs{E(G)} - 13 > \frac{10(n - s)}{7} - 13 > \frac{10(n' - s)}{7}.$$
\section{Proof of \texorpdfstring{\cref{lem:554}}{Lemma 2.1}}
\label{sec:554}
In this section we show our lemma on the case $D(e) \geq \langle 5,5,4\rangle$. Our proof follows similar techniques as in \cite{C21}, \cite{F21} and \cite{G21}. In particular, \cite{C21} analyzed the case $D(e) \geq \langle 4,4,4\rangle$, \cite{F21} analyzed the case $D(e) \geq \langle 5,5,5 \rangle$, and \cite{G21} analyzed the case $D(e) \geq \langle 5,5,3\rangle$. We use a slight variation of their methods to prove our lemma.

%%% TODO
Without loss of generality, assume $d(y), d(z) \geq 5$ and $d(x) \geq 4$. As we must not have $D(e) \geq \langle 6,4,2\rangle$, we must have $d(y) = d(z) = 5$. For $p \in \{x,y,z\}$, let $G(p)$ be the set of all vertices distinct from $x,y,z$ that lie on the same edge with $p$. We first note that we must have $G(y) = G(z)$. Suppose the contrary, and some edge $e_1 \neq e$ adjacent to $y$ contain some vertex not in $G(z)$. Then at most one edge adjacent to $z$ other than $e$ contains a vertex in $e_1$, so at least three edges $F$ adjacent to $z$ are disjoint from $e_1$. Thus, we can take an edge $e_2$ containing $x$ that is disjoint from $e_1$, then take an edge $e_3$ from $F$ that is disjoint from $e_2$. So $e,e_1,e_2,e_3$ forms a $C_{13}$, contradiction.

Similarly, we must have $G(x) \subset G(y)$. Suppose the contrary, and some edge $e_1 \neq e$ adjacent to $x$ contain some vertex not in $G(y)$. Then, we can take an edge $e_3$ containing $z$ that is disjoint from $e_1$. Among the four edges adjacent to $y$ distinct from $e$, at most two can intersect $e_3$, and at most one can intersect $e_1$. Thus, we can choose $e_2$ containing $y$ that is disjoint from $e_1$ and $e_3$. So $e,e_1,e_2,e_3$ forms a $C_{13}$, contradiction.

Thus $S \backslash \{x,y,z\} = G(y) = G(z) \supset G(x)$. We define $F$ as the set of all edges in $E(G)$ that contains one of the vertices in $S$, but is disjoint from $\{x,y,z\}$. It suffices to show that $F$ must be empty.

We denote the vertices in $G(z)$ by $a,b,c,d,r,s,p,q$, such that $\{z,a,b\}, \{z,c,d\}, \{z,r,s\}, \{z,p,q\}$ are edges in $E(G)$. 

\noindent {\it Step\uppercase\expandafter{\romannumeral 1}.}
We construct an auxiliary bipartite graph $H = (X_H,Y_H,E_H)$, where $X_H = \{e_i|y \in e_i\}, Y_H = \{e_j|z \in e_j\}, E_H = \{\{e_i,e_j\}| e_i \cap e_j \neq \emptyset\}$. $H$ is a 2-regular bipartite graph with order 8. Thus, $H=C_8$ or $H=C_4\biguplus C_4$.

We claim that if $G$ contains no crown, $H$ contains a $K_{2,2}$. Arbitrarily choose $e \in G(x)$. Define $V_1=e \cap S \subset E_H , W_1=\{e_i|e_i\cap V_1 \neq \emptyset\} \subset X_H \biguplus Y_H$, we have $|V_1|\le2, |W_1|\le4$, $|H-W_1|\ge 4$. To find a crown, we only need to choose $e_i \in X_G$ , $e_j \in Y_G$ s.t. $\{e_i, e_j\} \not\in  E_{G-W_1}$ . Therefore, if there is no crown in H,  $H-W_1$ has to be a completed bipartite graph. Since $|G-W_1|\ge4$ and two parts have the same order, there is definitely a $K_{2,2}$ in $H-W_1$. So $H$ contains a $K_{2,2}$, furthermore, $H=C_4\biguplus C_4$. 

By symmetry we can assume $\{z,a,b\}, \{z,c,d\}$ are in a $C_4$ and $\{z,r,s\}, \{z,p,q\}$ are in the other one. Without loss of generality we can further assume $\{y,b,d\}$, $\{y,a,c\}$ lie in $E(G)$, and $\{y,s,q\}, \{y,r,p\}$ lie in $E(G)$.

\noindent {\it Step\uppercase\expandafter{\romannumeral 2}.}
Now let $V_1 = \{a,b,c,d\}$, $V_2 = \{r,s,p,q\}$, We have symmetry between $V_1$ and $V_2$, and symmetry inside $V_i, i=1,2$ as well. We claim that there exists no edge containing $x$ that contains exactly one vertex in $V_1$ and another one in $V_2$. Otherwise we can let it be $\{x,a,r\}$ by symmetry. Then $\{z,a,b\}, \{y,b,d\}, \{z,p,q\}, \{x,a,r\}$ form a $C_{13}$, contradiction. Thus the edges other than $e$ containing $x$ must be a subset of $\{\{x,a,d\},\{x,b,c\},\{s,r,q\},\{x,s,p\}\}$.

\noindent {\it Step\uppercase\expandafter{\romannumeral 3}.}
Let $f$ be any element of $F$. By symmetry we can let $a \in f$ without loss of generality. Then we can see $b,c \notin f$. Firstly, we claim that $f$ cannot contain exactly one element $a$ of $S$. Otherwise $\{z,a,b\}, \{y,b,d\}, \{z,r,s\}, f$ form a $C_{13}$, contradiction. Secondly, we claim that $d \notin f$. Otherwise $G(x) = \{\{x,b,c\},\{s,r,q\},\{x,s,p\}\}$ since $d(x)\ge4$. Since at most one edge of $\{z,r,s\}$ and $\{z,p,q\}$ intersect $f$, we can assume $\{z,r,s\}\cap f = \emptyset$. Then$\{z,a,b\}, \{x,b,c\}, \{z,r,s\}, f$ form a $C_{13}$, contradiction.

Therefore we can assume $r\in f$ by symmetry. Similarly we know that $q \notin f$ since $a,d$ and $r,q$ are symmetric. So $f$ has exactly two elements $a,r$ of $S$. While $\{z,a,b\}, \{x,b,c\}, \{z,p,q\}, f$ form a $C_{13}$ in this case, contradiction.

\section*{Acknowledgements}
The main theorem is simultaneously obtained by Zhang and by Tang, Wu, and Zheng.

The second author's research is supported in part by National Natural Science Foundation of China grant 11931006, National Key Research and Development Program of China (Grant No. 2020YFA0713200), and the Shanghai Dawn Scholar Program grant 19SG01.

The third author's research is self-funded. 

The authors thank Professor Andr\'as Gy\'arf\'as for many valuable suggestions.

\end{document}